\documentclass[11pt]{article}
\usepackage{epsfig}

\def\be{\begin{equation}}
\def\ee{\end{equation}}
\def\bea{\begin{eqnarray}}
\def\eea{\end{eqnarray}}
\def\bes{\begin{eqnarray*}}
\def\ees{\end{eqnarray*}}

\def\nn{\nonumber}
\def\lb{\label}
\def\bs{\setminus}

\def\T{{\cal T}}
\def\H{{\cal H}}
\def\R{{\bf R}}
\def\C{{\bf C}}
\def\Z{{\bf Z}}

\def\N{{\bf N}}
\def\U{{\bf U}}

\def\Q{{\bf Q}}
\def\T{{\bf T}}

\def\Sg{{\Sigma}}

\def\aa{{\alpha}}
\def\bb{{\beta}}
\def\ga{{\gamma}}

\def\th{{\theta}}

\def\om{{\omega}}
\def\Om{{\Omega}}
\def\ep{{\epsilon}}
\def\lm{{\lambda}}
\def\Lm{{\Lambda}}

\def\sg{{\sigma}}
\def\dm{{\diamond}}
\def\Sg{{\Sigma}}
\def\vf{{\varphi}}

\def\<{{\langle}}
\def\>{{\rangle}}
\def\T{{\cal T}}

\def\P{{\cal P}}

\def\Nn{{\cal N}}

\def\mul{{\rm mul}}

\def\crit{{\rm crit}}

\def\Sp{{\rm Sp}}

\def\mod{{\rm mod}}

\def\wtd#1{\widetilde{#1}}

\def\hb{\vrule height0.18cm width0.14cm $\,$}

\title{A dichotomy result for closed characteristics on compact star-shaped hypersurfaces in $\R^{2n}$}

\author{Huagui Duan $^{1}$,\thanks{Partially supported by National Key R\&D Program of China (Grant No. 2020YFA0713300)
and NSFC (Nos. 11671215 and 11790271), LPMC of MOE of China and Nankai University. E-mail: duanhg@nankai.edu.cn.}
\quad Hui Liu $^{2}$,\thanks{Partially supported by NSFC (Nos. 12022111, 11771341). E-mail: huiliu00031514@whu.edu.cn.}
\quad Wenyan Ren $^{3}$ \thanks{E-mail: 2120210041@mail.nankai.edu.cn.}
\\\\
$^{1,3}$ School of Mathematical Sciences and LPMC, Nankai University, Tianjin 300071\\
$^{2}$ School of Mathematics and Statistics, Wuhan University,
Wuhan 430072, Hubei\\
The People's Republic of China \\}

\begin{document}
\date{}
\maketitle

\begin{abstract}
{\it In this paper, we prove that if all closed characteristics on a compact non-degenerate star-shaped hypersurface
$\Sigma$ in $\R^{2n}$ are elliptic, then either there exist exactly $n$  geometrically
distinct closed characteristics, or there exist infinitely many  geometrically
distinct closed characteristics.}
\end{abstract}

{\bf Key words}: Closed characteristic, star-shaped hypersurface, elliptic, Maslov-type index.

{\bf 2010 Mathematics Subject Classification}: 58E05, 37J45, 34C25.

\renewcommand{\theequation}{\thesection.\arabic{equation}}
\renewcommand{\thefigure}{\thesection.\arabic{figure}}

\setcounter{figure}{0}
\setcounter{equation}{0}
\section{Introduction and main results}

Let $\Sigma$ be a $C^3$ compact hypersurface in $\R^{2n}$ strictly star-shaped with respect to the origin, i.e.,
the tangent hyperplane at any $x\in\Sigma$ does not intersect the origin. We denote the set of all such
hypersurfaces by $\H_{st}(2n)$, and denote by $\H_{con}(2n)$ the subset of $\H_{st}(2n)$ which consists of all
strictly convex hypersurfaces. We consider closed characteristics $(\tau, y)$ on $\Sigma$, which are solutions
of the following problem
\be
\left\{\matrix{\dot{y}=JN_\Sigma(y), \cr
               y(\tau)=y(0), \cr }\right. \lb{1.1}\ee
where $J=\left(\matrix{0 &-I_n\cr
        I_n  & 0\cr}\right)$, $I_n$ is the identity matrix in $\R^n$, $\tau>0$, $N_\Sigma(y)$ is the outward
normal vector of $\Sigma$ at $y$ normalized by the condition $N_\Sigma(y)\cdot y=1$. Here $a\cdot b$ denotes
the standard inner product of $a, b\in\R^{2n}$. A closed characteristic $(\tau, y)$ is {\it prime}, if $\tau$
is the minimal period of $y$. Two closed characteristics $(\tau, y)$ and $(\sigma, z)$ are {\it geometrically
distinct}, if $y(\R)\not= z(\R)$. We denote by $\T(\Sigma)$ the set of geometrically distinct
closed characteristics $(\tau, y)$ on $\Sigma\in\mathcal{H}_{st}(2n)$. A closed characteristic
$(\tau,y)$ is {\it non-degenerate} if $1$ is a Floquet multiplier of $y$ of precisely algebraic multiplicity
$2$; {\it hyperbolic} if $1$ is a double Floquet multiplier of it and all the other Floquet multipliers
are not on ${\bf U}=\{z\in {\bf C}\mid |z|=1\}$, i.e., the unit circle in the complex plane; {\it elliptic}
if all the Floquet multipliers of $y$ are on ${\bf U}$; {\it irrationally elliptic} if $1$ is its double Floquet multipliers, and the other $(2n-2)$ locating on the unit circle with rotation angles being irrational multiples of $\pi$. We call a $\Sigma\in \mathcal{H}_{st}(2n)$ {\it non-degenerate} if
all the closed characteristics on $\Sigma$ together with all of their iterations are non-degenerate.

There is a long standing conjecture on the
number of closed characteristics on compact convex hypersurfaces in
$\R^{2n}$: \bea \,^{\#}\T(\Sg)\ge n, \qquad \forall \; \Sg\in\H_{con}(2n).
\lb{1.2}\eea

In 1978, P. Rabinowitz in \cite{Rab1} proved
$^\#\T(\Sg)\ge 1$ for any $\Sg\in\H_{st}(2n)$ and A. Weinstein in \cite{Wei1}
proved $^\#\T(\Sg)\ge 1$ for any $\Sg\in\H_{con}(2n)$ independently.
When $n\ge 2$,  in 1987-1988, I. Ekeland-L. Lassoued, I. Ekeland-H. Hofer, and A. Szulkin (cf. \cite{EkL1},
\cite{EkH1}, \cite{Szu1}) proved
$$ \,^{\#}\T(\Sg)\ge 2, \qquad \forall\,\Sg\in\H_{con}(2n). $$
In  \cite{HWZ1} of 1998, H. Hofer, K. Wysocki, and E. Zehnder proved $\,^{\#}\T(\Sg)=2$ or $\infty$ holds
for every $\Sg\in\H_{con}(4)$.
In \cite{LoZ} of 2002, Y. Long and C. Zhu further proved
\bea\;^{\#}\T(\Sg)\ge\left[\frac{n}{2}\right]+1, \qquad \forall\, \Sg\in \H_{con}(2n). \nn\eea
In particular, if all the prime closed characteristics on $\Sg$ are non-degenerate,
then $\,^{\#}\T(\Sg)\ge n$, cf. Theorem 1.1 and Corollary 1.1 of \cite{LoZ}.
 In \cite{WHL} of 2007,
W. Wang, X. Hu and Y. Long proved $\,^{\#}\T(\Sg)\ge 3$ for
every $\Sg\in\H_{con}(6)$.  In \cite{Wan2} of 2016, W. Wang proved $\,^{\#}\T(\Sg)\ge \left[\frac{n+1}{2}\right]+1$ for every $\Sg\in\H_{con}(2n)$. In \cite{Wan3} of 2016, W. Wang proved $\,^{\#}\T(\Sg)\ge 4$ for every $\Sg\in\H_{con}(8)$.

Note that every contact form supporting the standard contact structure on $M=S^{2n-1}$
arises from embeddings of $M$ into $\R^{2n}$ as a strictly star-shaped hypersurface enclosing the origin,
it is conjectured that in fact the conjecture (\ref{1.2}) holds for any $\Sg\in\H_{st}(2n)$. For star-shaped case, \cite{Gir1} of 1984 and \cite{BLMR} of 1985 show that $\;^{\#}\T(\Sg)\ge n$ for $\Sg\in\H_{st}(2n)$ under some pinching conditions. In \cite{HuL} of 2002, X. Hu and Y. Long proved that $\;^{\#}\T(\Sg)\ge 2$ for any non-degenerate $\Sg\in \H_{st}(2n)$. In \cite{HWZ2} of 2003,
H. Hofer, K. Wysocki, and E. Zehnder proved that $\,^{\#}\T(\Sg)=2$ or $\infty$ holds for every non-degenerate
$\Sg\in\H_{st}(4)$ provided that all stable and unstable manifolds of the hyperbolic closed orbits
on $\Sg$ intersect transversally, and recently this transversal condition was removed in an important paper by D. Cristofaro-Gardiner, M. Hutchings and D. Pomerleano \cite{CGHP}.  In \cite{CGH1} of 2016, D. Cristofaro-Gardiner and M. Hutchings proved that $\;^{\#}\T(\Sg)\ge 2$ for every contact manifold $\Sg$ of dimension three. Later various proofs of this result for star-shaped hypersurfaces have been given in \cite{GHHM}, \cite{LLo1} and \cite{GiG1}.

Since the appearance of \cite{HWZ1} and \cite{HWZ2}, Hofer, among others, has popularized in many talks
the following much more stronger conjecture than (\ref{1.2}), cf., Conjeture 1.5 of \cite{WHL} and Conjecture 1.1 of \cite{LLo3}:

{\bf Conjecture 1.1.} {\it For every integer $n\geq2$, there holds}
\bea \{\,^{\#}\T(\Sg)\mid \Sg\in\H_{st}(2n)\}=\{n\}\cup\{+\infty\}.\nn\eea

A typical example is the following weakly non-resonant ellipsoid: \bea \mathcal{E}_n(r)=\left\{z=(x_1, \ldots,x_n, y_1,\ldots,y_n)\in\R^{2n}\left |\frac{}{}\right.\frac{1}{2}\sum_{i=1}^n\frac{x_i^2+y_i^2}{r_i^2}=1\right\}\nn,\eea
where $r_i>0$ for $1\le i\le n$, we
have $\mathcal{E}_n(r)\in \mathcal{H}_{st}(2n)$. In this case, the corresponding Hamiltonian system is linear and all the solutions of (\ref{1.1}) can be computed explicitly. It is easy to verify that $^{\#}\T(\mathcal{E}_n(r)) = n$, and all the closed characteristics on $\mathcal{E}_n(r)$ are irrationally elliptic whenever $\frac{r_i}{r_j}\notin\Q$ for all $i\neq j$, and $^{\#}\T(\mathcal{E}_n(r)) = \infty$ otherwise.

Motivated by the dichotomy results of \cite{HWZ1} and \cite{HWZ2}, the detailed information about the ``exceptional'' case of hypersurfaces with exactly two closed characteristics arouses many interest. In \cite{Lon2} of 2000, Long proved that $\Sg\in\H_{con}(4)$
and $\,^{\#}\T(\Sg)=2$ imply that both of the closed characteristics must be elliptic.
In \cite{WHL} of 2007, W. Wang, X. Hu and Y. Long proved further that
$\Sg\in\H_{con}(4)$ and $\,^{\#}\T(\Sg)=2$ imply that both of the closed characteristics must be irrationally
elliptic. In \cite{LLo2} of 2015 and \cite{LLo3} of 2017, H. Liu and Y. Long proved
that the existence of exactly two closed characteristics on $\Sg\in\H_{st}(4)$ implies that both of them must be irrationally elliptic
provided that $\Sg$ is symmetric with respect to the origin. Recently, D. Cristofaro-Gardiner, U.L. Hryniewicz, M. Hutchings and H. Liu \cite{CGHHL}
removed this symmetric condition and proved a more general result that if there exist exactly two closed orbits on a contact manifold of dimension three, then both orbits are irrationally elliptic and this manifold is diffeomorphic to the three-sphere or a lens
space.

In \cite{Vit2} of 1989, C. Viterbo proved the existence of infinitely many closed characteristics on a compact star-shaped hypersurface $\Sigma$ in $\R^{4n}$, if all closed characteristics on $\Sigma$ are hyperbolic. This result has also been obtained by U. L. Hryniewicz and L. Macarini in Corollary 1.8 of \cite{HM}.

In this paper, motivated by the above results, we prove Conjecture 1.1 provided that all closed characteristics on $\Sigma\in\mathcal{H}_{st}(2n)$ are non-degenerate and elliptic.

\medskip

{\bf Theorem 1.2.} {\it If all closed characteristics on a compact non-degenerate star-shaped hypersurface
$\Sigma$ in $\R^{2n}$ are elliptic, then either there exist exactly $n$ geometrically
distinct closed characteristics, or there exist infinitely many geometrically
distinct closed characteristics.}

\medskip

In this paper, let $\N$, $\N_0$, $\Z$, $\Q$, $\R$, $\C$ and $\R^+$ denote the sets of natural integers,
non-negative integers, integers, rational numbers, real numbers, complex numbers and positive real
numbers respectively. We define the function $[a]=\max{\{k\in {\bf Z}\mid k\leq a\}}$, $\{a\}=a-[a]$,
$E(a)=\min{\{k\in{\bf Z}\mid k\geq a\}}$ and $\varphi(a)=E(a)-[a]$. Denote by $a\cdot b$ and $|a|$ the standard
inner product and norm in $\R^{2n}$. Denote by $\langle\cdot,\cdot\rangle$ and $\|\cdot\|$
the standard $L^2$ inner product and $L^2$ norm. For an $S^1$-space $X$, we denote by
$X_{S^1}$ the homotopy quotient of $X$ by $S^1$, i.e., $X_{S^1}=S^\infty\times_{S^1}X$,
where $S^\infty$ is the unit sphere in an infinite dimensional {\it complex} Hilbert space.
In this paper we use $\Q$ coefficients for all homological and cohomological modules. By $t\to a^+$, we
mean $t>a$ and $t\to a$.

\setcounter{figure}{0}
\setcounter{equation}{0}
\section{Mean index identities for closed characteristics on compact star-shaped hypersurfaces in $\R^{2n}$}

In this section, we briefly review the mean index identities for
closed characteristics on $\Sg\in\H_{st}(2n)$ developed in \cite{LLW14} which will be needed in Section 4. All
the details of proofs can be found in \cite{LLW14}.
Now we fix a $\Sg\in\H_{st}(2n)$ and assume the following condition on $\T(\Sg)$:

\medskip

(F) {\bf There exist only finitely many geometrically distinct prime closed characteristics\\
$\qquad\qquad \{(\tau_j, y_j)\}_{1\le j\le k}$ on $\Sigma$. }

\medskip

Let $\hat{\tau}=\inf_{1\leq j\leq k}{\tau_j}$ and $T$ be a fixed positive constant. Then by Section 2 of
\cite{LLW14}, for any $a>\frac{\hat{\tau}}{T}$, we can construct a function $\varphi_a\in C^{\infty}({\bf R}, {\bf R}^+)$
which has $0$ as its unique critical point in $[0, +\infty)$. Moreover, $\frac{\varphi^{\prime}(t)}{t}$ is strictly
decreasing for $t>0$ together with $\varphi(0)=0=\varphi^{\prime}(0)$ and
$\varphi^{\prime\prime}(0)=1=\lim_{t\rightarrow 0^+}\frac{\varphi^{\prime}(t)}{t}$. More precisely, we
define $\varphi_a$ and the Hamiltonian function $\wtd{H}_a(x)=a\vf_a(j(x))$ via Lemma 2.2 and Lemma 2.4
in \cite{LLW14}. The precise dependence of $\varphi_a$ on $a$ is explained in Remark 2.3 of \cite{LLW14}.

For technical reasons we further modify the Hamiltonian $\wtd{H}_a$, and define the new Hamiltonian
function $H_a$ via Proposition 2.5 of \cite{LLW14} and consider the fixed period problem
\be  \dot{x}(t)=JH_a^\prime(x(t)),\quad x(0)=x(T).  \lb{2.1}\ee
Then $H_a\in C^{3}({\bf R}^{2n} \setminus\{0\},{\bf R})\cap C^{1}({\bf R}^{2n},{\bf R})$.
Solutions of (\ref{2.1}) are $x\equiv 0$ and $x=\rho y(\tau t/T)$ with
$\frac{\vf_a^\prime(\rho)}{\rho}=\frac{\tau}{aT}$, where $(\tau, y)$ is a solution of (\ref{1.1}). In particular,
non-zero solutions of (\ref{2.1}) are in one to one correspondence with solutions of (\ref{1.1}) with period
$\tau<aT$.

For any $a>\frac{\hat{\tau}}{T}$, we can choose some large constant $K=K(a)$ such that
\be H_{a,K}(x) = H_a(x)+\frac{1}{2}K|x|^2   \lb{2.2}\ee
is a strictly convex function, that is,
\be (\nabla H_{a, K}(x)-\nabla H_{a, K}(y), x-y) \geq \frac{\ep}{2}|x-y|^2,  \lb{2.3}\ee
for all $x, y\in {\bf R}^{2n}$, and some positive $\ep$. Let $H_{a,K}^*$ be the Fenchel dual of $H_{a,K}$
defined by
$$  H_{a,K}^\ast (y) = \sup\{x\cdot y-H_{a,K}(x)\;|\; x\in \R^{2n}\}.   $$
The dual action functional on $X=W^{1, 2}({\bf R}/{T {\bf Z}}, {\bf R}^{2n})$ is defined by
\be F_{a,K}(x) = \int_0^T{\left[\frac{1}{2}(J\dot{x}-K x,x)+H_{a,K}^*(-J\dot{x}+K x)\right]dt}. \lb{2.4}\ee
Then $F_{a,K}\in C^{1,1}(X, \R)$ and for $KT\not\in 2\pi{\bf Z}$, $F_{a,K}$ satisfies the
Palais-Smale condition and $x$ is a critical point of $F_{a, K}$ if and only if it is a solution of (\ref{2.1}). Moreover,
$F_{a, K}(x_a)<0$ and it is independent of $K$ for every critical point $x_a\neq 0$ of $F_{a, K}$.

When $KT\notin 2\pi{\bf Z}$, the map $x\mapsto -J\dot{x}+Kx$ is a Hilbert space isomorphism between
$X=W^{1, 2}({\bf R}/({T {\bf Z}}); {\bf R}^{2n})$ and $E=L^{2}({\bf R}/(T {\bf Z}),{\bf R}^{2n})$. We denote its inverse
by $M_K$ and the functional
\be \Psi_{a,K}(u)=\int_0^T{\left[-\frac{1}{2}(M_{K}u, u)+H_{a,K}^*(u)\right]dt}, \qquad \forall\,u\in E. \lb{2.5}\ee
Then $x\in X$ is a critical point of $F_{a,K}$ if and only if $u=-J\dot{x}+Kx$ is a critical point of $\Psi_{a, K}$.

Suppose $u$ is a nonzero critical point of $\Psi_{a, K}$.
Then the formal Hessian of $\Psi_{a, K}$ at $u$ is defined by
\be Q_{a,K}(v)=\int_0^T(-M_K v\cdot v+H_{a,K}^{*\prime\prime}(u)v\cdot v)dt,  \lb{2.6}\ee
which defines an orthogonal splitting $E=E_-\oplus E_0\oplus E_+$ of $E$ into negative, zero and positive subspaces.
The index and nullity of $u$ are defined by $i_K(u)=\dim E_-$ and $\nu_K(u)=\dim E_0$ respectively.
Similarly, we define the index and nullity of $x=M_Ku$ for $F_{a, K}$, we denote them by $i_K(x)$ and
$\nu_K(x)$. Then we have
\be  i_K(u)=i_K(x),\quad \nu_K(u)=\nu_K(x),  \lb{2.7}\ee
which follow from the definitions (\ref{2.4}) and (\ref{2.5}). The following important formula was proved in
Lemma 6.4 of \cite{Vit2}:
\be  i_K(x) = 2n([KT/{2\pi}]+1)+i^v(x) \equiv d(K)+i^v(x),   \lb{2.8}\ee
where the Viterbo index $i^v(x)$ does not depend on K, but only on $H_a$.

By the proof of Proposition 2 of \cite{Vit1}, we have that $v\in E$ belongs to the null space of $Q_{a, K}$
if and only if $z=M_K v$ is a solution of the linearized system
\be  \dot{z}(t) = JH_a''(x(t))z(t).  \lb{2.9}\ee
Thus the nullity in (\ref{2.7}) is independent of $K$, which we denote by $\nu^v(x)\equiv \nu_K(u)= \nu_K(x)$.

By Proposition 2.11 of \cite{LLW14}, the index $i^v(x)$ and nullity $\nu^v(x)$ coincide with those defined for
the Hamiltonian $H(x)=j(x)^\alpha$ for all $x\in\R^{2n}$ and some $\aa\in (1,2)$. Especially
$1\le \nu^v(x)\le 2n-1$ always holds.

For every closed characteristic $(\tau, y)$ on $\Sigma$, let $aT>\tau$ and choose $\vf_a$ as above.
Determine $\rho$ uniquely by $\frac{\vf_a'(\rho)}{\rho}=\frac{\tau}{aT}$. Let $x=\rho y(\frac{\tau t}{T})$.
Then we define the index $i(\tau,y)$ and nullity $\nu(\tau,y)$ of $(\tau,y)$ by
$$ i(\tau,y)=i^v(x), \qquad \nu(\tau,y)=\nu^v(x). $$
Then the mean index of $(\tau,y)$ is defined by
\bea \hat i(\tau,y) = \lim_{m\rightarrow\infty}\frac{i(m\tau,y)}{m}.  \nn\eea
Note that by Proposition 2.11 of \cite{LLW14}, the index and nullity are well defined and are independent of the
choice of $a$. For a closed characteristic $(\tau,y)$ on $\Sigma$, we simply denote by $y^m\equiv(m\tau,y)$
the m-th iteration of $y$ for $m\in\N$.

We have a natural $S^1$-action on $X$ or $E$ defined by
$$  \theta\cdot u(t)=u(\theta+t),\quad\forall\, \theta\in S^1, \, t\in\R.  $$
Clearly both of $F_{a, K}$ and $\Psi_{a, K}$ are $S^1$-invariant. For any $\kappa\in\R$, we denote by
\bea
\Lambda_{a, K}^\kappa &=& \{u\in L^{2}({\bf R}/({T {\bf Z}}); {\bf R}^{2n})\;|\;\Psi_{a,K}(u)\le\kappa\}  \nn\\
X_{a, K}^\kappa &=& \{x\in W^{1, 2}({\bf R}/(T {\bf Z}),{\bf R}^{2n})\;|\;F_{a, K}(x)\le\kappa\}.  \nn\eea
For a critical point $u$ of $\Psi_{a, K}$ and the corresponding $x=M_K u$ of $F_{a, K}$, let
\bea
\Lm_{a,K}(u) &=& \Lm_{a,K}^{\Psi_{a, K}(u)}
   = \{w\in L^{2}(\R/(T\Z), \R^{2n}) \;|\; \Psi_{a, K}(w)\le\Psi_{a,K}(u)\},  \nn\\
X_{a,K}(x) &=& X_{a,K}^{F_{a,K}(x)} = \{y\in W^{1, 2}(\R/(T\Z), \R^{2n}) \;|\; F_{a,K}(y)\le F_{a,K}(x)\}. \nn\eea
Clearly, both sets are $S^1$-invariant. Denote by $\crit(\Psi_{a, K})$ the set of critical points of $\Psi_{a, K}$.
Because $\Psi_{a,K}$ is $S^1$-invariant, $S^1\cdot u$ becomes a critical orbit if $u\in \crit(\Psi_{a, K})$.
Note that by the condition (F), the number of critical orbits of $\Psi_{a, K}$
is finite. Hence as usual we can make the following definition.

{\bf Definition 2.1.} {\it Suppose $u$ is a nonzero critical point of $\Psi_{a, K}$, and $\Nn$ is an $S^1$-invariant
open neighborhood of $S^1\cdot u$ such that $\crit(\Psi_{a,K})\cap (\Lm_{a,K}(u)\cap \Nn) = S^1\cdot u$.
Then the $S^1$-critical module of $S^1\cdot u$ is defined by
$$ C_{S^1,\; q}(\Psi_{a, K}, \;S^1\cdot u)
=H_{q}((\Lambda_{a, K}(u)\cap\Nn)_{S^1},\; ((\Lambda_{a,K}(u)\setminus S^1\cdot u)\cap\Nn)_{S^1}). $$
Similarly, we define the $S^1$-critical module $C_{S^1,\; q}(F_{a, K}, \;S^1\cdot x)$ of $S^1\cdot x$
for $F_{a, K}$.}

We fix $a$ and let $u_K\neq 0$ be a critical point of $\Psi_{a, K}$ with multiplicity $\mul(u_K)=m$,
that is, $u_K$ corresponds to a closed characteristic $(\tau, y)\subset\Sigma$ with $(\tau, y)$
being $m$-iteration of
some prime closed characteristic. Precisely, we have $u_K=-J\dot x+Kx$ with $x$
being a solution of (\ref{2.1}) and $x=\rho y(\frac{\tau t}{T})$ with
$\frac{\vf_a^\prime(\rho)}{\rho}=\frac{\tau}{aT}$.
Moreover, $(\tau, y)$ is a closed characteristic on $\Sigma$ with minimal period $\frac{\tau}{m}$.
For any $p\in\N$ satisfying $p\tau<aT$, we choose $K$
such that $pK\notin \frac{2\pi}{T}\Z$, then the $p$th iteration $u_{pK}^p$ of $u_K$ is given by $-J\dot x^p+pKx^p$,
where $x^p$ is the unique solution of (\ref{2.1}) corresponding to $(p\tau, y)$
and is a critical point of $F_{a, pK}$, that
is, $u_{pK}^p$ is the critical point of $\Psi_{a, pK}$ corresponding to $x^p$.

\medskip

{\bf Lemma 2.2.} (cf. Proposition 4.2 and Remark 4.4 of \cite{LLW14} ) {\it If $u_{pK}^p$ is non-degenerate,
i.e., $\nu_{pK}(u_{pK}^p)=1$, let
$\bb(x^p)=(-1)^{i_{pK}(u_{pK}^p)-i_{K}(u_{K})}=(-1)^{i^v(x^p)-i^v(x)}$, then}
\bea C_{S^1,q-d(pK)+d(K)}(F_{a,K},S^1\cdot x^p)
&=& C_{S^1,q}(F_{a,pK},S^1\cdot x^p)=C_{S^1,q}(\Psi_{a,pK},S^1\cdot u^p_{pK}) \nn\\
&=& \left\{\matrix{
     \Q, &\quad {\it if}\;\; q=i_{pK}(u_{pK}^p),\;\;{\it and}\;\;\bb(x^p)=1, \cr
     0, &\quad {\it otherwise}. \cr}\right.  \lb{2.10}\eea

\medskip

{\bf Theorem 2.3.} (cf. Theorem 1.1 of \cite{LLW14} and Theorem 1.2 of \cite{Vit2}) {\it Suppose that
$\Sg\in\H_{st}(2n)$ satisfying $^\#\T(\Sg)<+\infty$. Denote by $\{(\tau_j,y_j)\}_{1\le j\le k}$ all the
geometrically distinct prime closed characteristics. Then the following identities hold
\bea \sum_{1\le j\le k \atop \hat{i}(y_j)>0}\frac{\hat{\chi}(y_j)}{\hat{i}(y_j)}=\frac{1}{2},\qquad
\sum_{1\le j\le k \atop \hat{i}(y_j)<0}\frac{\hat{\chi}(y_j)}{\hat{i}(y_j)}=0,\lb{2.11}\eea
where $\hat{\chi}(y)\in\Q$ is the average Euler characteristic given by Definition 4.8 and Remark 4.9 of \cite{LLW14}.

In particular, if all $y^m$'s are non-degenerate for $m\ge 1$, then
\bea \hat{\chi}(y)=\left\{\matrix{
     (-1)^{i(y)}, &\quad {\it if}\;\; i(y^2)-i(y)\in 2\Z, \cr
     \frac{(-1)^{i(y)}}{2}, &\quad {\it otherwise}. \cr}\right.  \lb{2.12}\eea}

Let $F_{a, K}$ be a functional defined by (\ref{2.4}) for some $a, K\in\R$ large enough and let $\ep>0$ be
small enough such that $[-\ep, 0)$ contains no critical values of $F_{a, K}$. For $b$ large enough,
The normalized Morse series of $F_{a, K}$ in $ X^{-\ep}\setminus X^{-b}$
is defined, as usual, by
\be  M_a(t)=\sum_{q\ge 0,\;1\le j\le p} \dim C_{S^1,\;q}(F_{a, K}, \;S^1\cdot v_j)t^{q-d(K)},  \lb{2.13}\ee
where we denote by $\{S^1\cdot v_1, \ldots, S^1\cdot v_p\}$ the critical orbits of $F_{a, K}$ with critical
values less than $-\ep$. The Poincar\'e series of $H_{S^1, *}( X, X^{-\ep})$ is $t^{d(K)}Q_a(t)$, according
to Theorem 5.1 of \cite{LLW14}, if we set $Q_a(t)=\sum_{k\in \Z}{q_kt^k}$, then
$$   q_k=0, \qquad \forall\;k\in \mathring {I},  $$
where $I$ is an interval of $\Z$ such that $I \cap [i(\tau, y), i(\tau, y)+\nu(\tau, y)-1]=\emptyset$ for all
closed characteristics $(\tau,\, y)$ on $\Sigma$ with $\tau\ge aT$. Then by Section 6 of \cite{LLW14}, we have
$$  M_a(t)-\frac{1}{1-t^2}+Q_a(t) = (1+t)U_a(t),   $$
where $U_a(t)=\sum_{i\in \Z}{u_it^i}$ is a Laurent series with nonnegative coefficients.
If there is no closed characteristic with $\hat{i}=0$, then
\be   M(t)-\frac{1}{1-t^2}=(1+t)U(t),    \lb{2.14}\ee
where $M(t)=\sum_{p\in \Z}{M_pt^p}$ denotes $M_a(t)$ as $a$ tends to infinity. In addition, we also denote by $b_p$ the coefficient of $t^p$ of $\frac{1}{1-t^2}=\sum_{p\in \Z}{b_pt^p}$, i.e. there holds $b_p=1$, $\forall\ p\in2\N_0$ and $b_p=0$, $\forall\ p\not\in2\N_0$.

\setcounter{figure}{0}
\setcounter{equation}{0}
\section{The Maslov-type index theory of symplectic paths}

In \cite{Lon2} of 1999, Y. Long established the basic normal form
decomposition of symplectic matrices. Based on this result he
further established the precise iteration formulae of indices of
symplectic paths in \cite{Lon3} of 2000.

As in \cite{Lon3}, denote by
\bea
N_1(\lm, b) &=& \left(\matrix{\lm & b\cr
                                0 & \lm\cr}\right), \qquad {\rm for\;}\lm=\pm 1, \; b\in\R, \nn\\
D(\lm) &=& \left(\matrix{\lm & 0\cr
                      0 & \lm^{-1}\cr}\right), \qquad {\rm for\;}\lm\in\R\bs\{0, \pm 1\}, \nn\\
R(\th) &=& \left(\matrix{\cos\th & -\sin\th \cr
                           \sin\th & \cos\th\cr}\right), \qquad {\rm for\;}\th\in (0,\pi)\cup (\pi,2\pi), \nn\\
N_2(e^{\th\sqrt{-1}}, B) &=& \left(\matrix{ R(\th) & B \cr
                  0 & R(\th)\cr}\right), \qquad {\rm for\;}\th\in (0,\pi)\cup (\pi,2\pi)\;\; {\rm and}\; \nn\\
        && \qquad B=\left(\matrix{b_1 & b_2\cr
                                  b_3 & b_4\cr}\right)\; {\rm with}\; b_j\in\R, \;\;
                                         {\rm and}\;\; b_2\not= b_3. \nn\eea
Here $N_2(e^{\th\sqrt{-1}}, B)$ is non-trivial if $(b_2-b_3)\sin\theta<0$, and trivial
if $(b_2-b_3)\sin\theta>0$.

As in \cite{Lon3}, the $\diamond$-sum (direct sum) of any two real matrices is defined by
$$ \left(\matrix{A_1 & B_1\cr C_1 & D_1\cr}\right)_{2i\times 2i}\diamond
      \left(\matrix{A_2 & B_2\cr C_2 & D_2\cr}\right)_{2j\times 2j}
=\left(\matrix{A_1 & 0 & B_1 & 0 \cr
                                   0 & A_2 & 0& B_2\cr
                                   C_1 & 0 & D_1 & 0 \cr
                                   0 & C_2 & 0 & D_2}\right). $$

For every $M\in\Sp(2n)$, the homotopy set $\Omega(M)$ of $M$ in $\Sp(2n)$ is defined by
$$ \Om(M)=\{N\in\Sp(2n)\,|\,\sg(N)\cap\U=\sg(M)\cap\U\equiv\Gamma\;\mbox{and}
                    \;\nu_{\om}(N)=\nu_{\om}(M)\, \forall\om\in\Gamma\}, $$
where $\sg(M)$ denotes the spectrum of $M$,
$\nu_{\om}(M)\equiv\dim_{\C}\ker_{\C}(M-\om I)$ for $\om\in\U$.
The component $\Om^0(M)$ of $P$ in $\Sp(2n)$ is defined by
the path connected component of $\Om(M)$ containing $M$.

\medskip

{\bf Lemma 3.1.} (cf. \cite{Lon3}, Lemma 9.1.5 and List 9.1.12 of \cite{Lon4})
{\it For $M\in\Sp(2n)$ and $\om\in\U$, the splitting number $S_M^\pm(\om)$ (cf. Definition 9.1.4 of \cite{Lon4}) satisfies
\begin{eqnarray}
S_M^{\pm}(\om) &=& 0, \qquad {\it if}\;\;\om\not\in\sg(M).  \nn\\
S_{N_1(1,a)}^+(1) &=& \left\{\matrix{1, &\quad {\rm if}\;\; a\ge 0, \cr
0, &\quad {\rm if}\;\; a< 0. \cr}\right. \nn\eea

For any $M_i\in\Sp(2n_i)$ with $i=0$ and $1$, there holds }
\bea S^{\pm}_{M_0\diamond M_1}(\om) = S^{\pm}_{M_0}(\om) + S^{\pm}_{M_1}(\om),
    \qquad \forall\;\om\in\U. \nn\eea

\medskip

For every $\ga\in\mathcal{P}_\tau(2n)\equiv\{\ga\in C([0,\tau],Sp(2n))\ |\ \ga(0)=I_{2n}\}$, we extend
$\ga(t)$ to $t\in [0,m\tau]$ for every $m\in\N$ by
\bea \ga^m(t)=\ga(t-j\tau)\ga(\tau)^j \qquad \forall\;j\tau\le t\le (j+1)\tau \;\;
               {\rm and}\;\;j=0, 1, \ldots, m-1, \lb{3.1}\eea
as in P.114 of \cite{Lon2}. As in \cite{LoZ} and \cite{Lon4}, we denote the Maslov-type indices of
$\ga^m$ by $(i(\ga,m),\nu(\ga,m))$.

The following is the precise index iteration formulae for
symplectic paths, which is due to Y. Long
(cf. Chapter 8 of \cite{Lon4} or Theorems 6.5 and 6.7 of \cite{LoZ}).

\medskip

{\bf Theorem 3.2.} {\it Let $\ga\in\P_{\tau}(2n)$. Then there exists a path $f\in
C([0,1],\Omega^0(\gamma(\tau))$ such that $f(0)=\gamma(\tau)$ and
\bea f(1)&=&N_1(1,1)^{\diamond p_-} \diamond I_{2p_0}\diamond
N_1(1,-1)^{\diamond p_+}
\diamond N_1(-1,1)^{\diamond q_-} \diamond (-I_{2q_0})\diamond
N_1(-1,-1)^{\diamond q_+}\nn\\
&&\diamond R(\theta_1)\diamond\cdots\diamond R(\theta_r)
\diamond N_2(\omega_1, u_1)\diamond\cdots\diamond N_2(\omega_{r_*}, u_{r_*}) \nn\\
&&\diamond N_2(\lm_1, v_1)\diamond\cdots\diamond N_2(\lm_{r_0}, v_{r_0})
\diamond M_0 \nn\eea
where $ N_2(\omega_j, u_j) $s are
non-trivial and   $ N_2(\lm_j, v_j)$s  are trivial basic normal
forms; $\sigma (M_0)\cap U=\emptyset$; $p_-$, $p_0$, $p_+$, $q_-$,
$q_0$, $q_+$, $r$, $r_*$ and $r_0$ are non-negative integers;
$\omega_j=e^{\sqrt{-1}\alpha_j}$, $
\lambda_j=e^{\sqrt{-1}\beta_j}$; $\theta_j$, $\alpha_j$, $\beta_j$
$\in (0, \pi)\cup (\pi, 2\pi)$; these integers and real numbers
are uniquely determined by $\gamma(\tau)$. Then using the
functions defined in (\ref{1.2}), we have
\bea i(\gamma, m)&=&m(i(\gamma,
1)+p_-+p_0-r)+2\sum_{j=1}^r E\left(\frac{m\theta_j}{2\pi}\right)-r
-p_--p_0\nn\\&&-\frac{1+(-1)^m}{2}(q_0+q_+)+2\left(
\sum_{j=1}^{r_*}\varphi\left(\frac{m\alpha_j}{2\pi}\right)-r_*\right),
\nn\eea
\bea \nu(\gamma, m)&=&\nu(\gamma,
1)+\frac{1+(-1)^m}{2}(q_-+2q_0+q_+)+2(r+r_*+r_0)\nn\\
&&-2\left(\sum_{j=1}^{r}\varphi\left(\frac{m\theta_j}{2\pi}\right)+
\sum_{j=1}^{r_*}\varphi\left(\frac{m\alpha_j}{2\pi}\right)
+\sum_{j=1}^{r_0}\varphi\left(\frac{m\beta_j}{2\pi}\right)\right),\nn\\
\hat i(\gamma, 1)&=& i(\gamma, 1)+p_-+p_0-r+\sum_{j=1}^r
\frac{\theta_j}{\pi},\nn\eea
where $N_1(1, \pm 1)=
\left(\matrix{ 1 &\pm 1\cr 0 & 1\cr}\right)$, $N_1(-1, \pm 1)=
\left(\matrix{ -1 &\pm 1\cr 0 & -1\cr}\right)$,
$R(\theta)=\left(\matrix{\cos\th &
                  -\sin\th\cr\sin\th & \cos\th\cr}\right)$,
$ N_2(\omega, b)=\left(\matrix{R(\th) & b
                  \cr 0 & R(\th)\cr}\right)$    with some
$\th\in (0,\pi)\cup (\pi,2\pi)$ and $b=
\left(\matrix{ b_1 &b_2\cr b_3 & b_4\cr}\right)\in\R^{2\times2}$,
such that $(b_2-b_3)\sin\theta>0$, if $ N_2(\omega, b)$ is
trivial; $(b_2-b_3)\sin\theta<0$, if $ N_2(\omega, b)$ is
non-trivial. We have $i(\gamma, 1)$ is odd if $f(1)=N_1(1, 1)$, $I_2$,
$N_1(-1, 1)$, $-I_2$, $N_1(-1, -1)$ and $R(\theta)$; $i(\gamma, 1)$ is
even if $f(1)=N_1(1, -1)$ and $ N_2(\omega, b)$; $i(\gamma, 1)$ can be any
integer if $\sigma (f(1)) \cap \U=\emptyset$.}

\medskip

The common index jump theorem (cf. Theorem 4.3 of \cite{LoZ}) for symplectic paths established by Long
and Zhu in 2002 has become one of the main tools to study the multiplicity and stability problems of
closed solution orbits in Hamiltonian and symplectic dynamics. Recently, the following enhanced version
of it has been obtained by Duan, Long and Wang in \cite{DLW}, which will play an important role in the
proof in Section 4.

\medskip

{\bf Theorem 3.3.} (cf. Theorem 3.5 of \cite{DLW}) ({\bf The enhanced common index jump theorem for
symplectic paths}) {\it Let $\gamma_k\in\mathcal{P}_{\tau_k}(2n)$ for $k=1,\cdots,q$ be a finite
collection of symplectic paths. Let $M_k=\ga_k(\tau_k)$. We extend $\ga_k$ to $[0,+\infty)$ by (\ref{3.1})
inductively. Suppose
\bea  \hat{i}(\ga_k,1) > 0, \qquad \forall\ k=1,\cdots,q.  \nn\eea
Then for every integer $\bar{m}\in \N$, there exist infinitely many $(q+1)$-tuples
$(N, m_1,\cdots,m_q) \in \N^{q+1}$ such that for all $1\le k\le q$ and $1\le m\le \bar{m}$, there holds
\bea
\nu(\ga_k,2m_k-m) &=& \nu(\ga_k,2m_k+m) = \nu(\ga_k, m),   \nn\\
i(\ga_k,2m_k+m) &=& 2N+i(\ga_k,m),                         \nn\\
i(\ga_k,2m_k-m) &=& 2N-i(\ga_k,m)-2(S^+_{M_k}(1)+Q_k(m)),  \nn\\
i(\ga_k, 2m_k)&=& 2N -(S^+_{M_k}(1)+C(M_k)-2\Delta_k),     \nn\eea
where $C(M_k)=\sum_{0<\th<2\pi}S_{M_k}^-(e^{\sqrt{-1}\th})$ and}
\bea \Delta_k = \sum_{0<\{m_k\th/\pi\}<\delta}S^-_{M_k}(e^{\sqrt{-1}\th}),\qquad
 Q_k(m) = \sum_{e^{\sqrt{-1}\th}\in\sg(M_k),\atop \{\frac{m_k\th}{\pi}\}
                   = \{\frac{m\th}{2\pi}\}=0}S^-_{M_k}(e^{\sqrt{-1}\th}). \nn\eea

The following is the relation between Viterbo index and Maslov-type index of symplectic path.

{\bf Theorem 3.4.} (cf. Theorem 2.1 of \cite{HuL} and Theorem 6.1 of \cite{LLo2}) {\it Suppose $\Sg\in \H_{st}(2n)$ and
$(\tau,y)\in \T(\Sigma)$. Then we have
\bea i(y^m)\equiv i(m\tau,y)=i(y, m)-n,\quad \nu(y^m)\equiv\nu(m\tau, y)=\nu(y, m),
       \qquad \forall m\in\N, \nn\eea
where $i(y^m)$ and $\nu(y^m)$ are the index and nullity of $(m\tau,y)$ defined in Section 2, $i(y, m)$ and $\nu(y, m)$
are the Maslov-type index and nullity of $(m\tau,y)$ (cf. Section 5.4 of \cite{Lon3}). In particular, we have
$\hat{i}(\tau,y)=\hat{i}(y,1)$, where $\hat{i}(\tau ,y)$ is given in Section 2, $\hat{i}(y,1)$
is the mean Maslov-type index (cf. Definition 8.1 of \cite{Lon4}). Hence we denote it simply by $\hat{i}(y)$.}

\setcounter{figure}{0}
\setcounter{equation}{0}
\section{Proof of Theorem 1.2}

In order to prove Theorem 1.2, we make the following assumption

\medskip

{\bf (ECC)} {\it Suppose that all prime closed characteristics on a compact non-degenerate star-shaped hypersurface $\Sigma$
in $\R^{2n}$ are elliptic, and the total number of distinct closed characteristics on $\Sg$ is finite, denoted by $\{(\tau_k,y_k)\}_{k=1}^q$.}

\medskip

We denote by $\ga_k\equiv \ga_{y_k}$ the associated symplectic path of $(\tau_k,y_k)$ for $1\le k\le q$ in the assumption (ECC).
Then by Lemma 3.3 of \cite{HuL},
there exists $P_k\in Sp(2n)$ and $U_k\in Sp(2n-2)$ such that
\bea M_k\equiv\ga_k(\tau_k)=P_k^{-1}(N_1(1,1)\dm U_k)P_k,\qquad 1\le k\le q,\lb{4.1a}\eea
where every $U_k$ has the following form:
\bea
U_k &=& R(\th^k_{1})\,\dm\,\cdots\,\dm\,R(\th^k_{r_k}) \nn\\
& & \dm\,N_2(e^{\aa^k_{1}\sqrt{-1}},A^k_{1})\,\dm\,\cdots\,\dm\,N_2(e^{\aa^k_{r_{k\ast}}\sqrt{-1}},A^k_{r_{k\ast}})\,   \dm\,N_2(e^{\bb^k_{1}\sqrt{-1}},B^k_{1})\,\dm\,\cdots\,\dm\,N_2(e^{\bb^k_{r_{k0}}\sqrt{-1}},B^k_{r_{k0}}), \nn\eea
where $\frac{\th^k_{j}}{2\pi}\in(0,1)\backslash\Q$ for $1\le j\le r_k$,
$\frac{\aa^k_{j}}{2\pi}\in (0,1)\backslash\Q$ for $1\le j\le r_{k\ast}$, $\frac{\bb^k_{j}}{2\pi}\in (0,1)\backslash\Q$ for $1\le j\le r_{k0}$,
and
\be r_k + 2r_{k\ast} + 2r_{k0} = n - 1. \lb{4.2a}\ee
By Lemma 3.1 and Theorem 3.2 we obtain the index iteration formula of $y_k^m$
\be i(y_k,m) = m(i(y_k,1)+1-r_k) + 2\sum_{j=1}^{r_k}\left[\frac{m\th^k_j}{2\pi}\right] +r_k-1,
    \quad\forall\  1\le k\le q,\ m\ge 1. \lb{4.3a}\ee
Therefore by Theorem 3.4 it yields
\be i(y_k^m) = m(i(y_k)+n+1-r_k) + 2\sum_{j=1}^{r_k}\left[\frac{m\th^k_j}{2\pi}\right] +r_k-n-1,
    \quad\forall\  1\le k\le q,\ m\ge 1, \lb{4.4a}\ee
which implies
\be \hat{i}(y_k,1)=\hat{i}(y_k) = i(y_k)+n+1-r_k + \sum_{j=1}^{r_k}\frac{\th^k_j}{\pi},
    \quad\forall\ 1\le k\le q, \lb{4.5a}\ee

\medskip

{\bf Claim 1.} {\it For any $1\le k\le q$ and $m\ge 1$, there holds $i(y_k^m)\in 2\Z$.}

\medskip

In fact, firstly by Theorem 3.2, every path $\ga\in\P_{\tau}(2)$ with its end matrix homotopic to $N_1(1,1)$ and $R(\theta)$ has odd index $i(\ga,1)$, and every path $\ga\in\P_{\tau}(4)$ with its end matrix homotopic to $N_2(\om,B)$ has even indices $i(\ga,1)$. Then by (\ref{4.1a}) and the homotopy invariance and the symplectic additivity of Maslov-type indices (cf. Theorem 6.2.7 of \cite{Lon4}), for every closed characteristic $y_k$, $1\le k\le q$, $i(y_k,1)$ has the same parity as $n$. Thus by Theorem 3.4, it yields
\bea i(y_k)=i(y_k,1)-n=0 \quad (\mod\,2), \qquad \forall\ 1\le k\le q. \lb{4.6a}\eea

By (\ref{4.2a}), (\ref{4.4a}) and (\ref{4.6a}), we get
\bea i(y_k^{m+1})-i(y_k^m)&=&(i(y_k)+n+1-r_k) + 2\sum_{j=1}^{r_k}\left[\frac{(m+1)\th^k_j}{2\pi}\right]-2\sum_{j=1}^{r_k}\left[\frac{m\th^k_j}{2\pi}\right] \nn\\
&=&i(y_k)+(n+1)+[2r_{k\ast} + 2r_{k0}-(n-1)]\quad (\mod\,2)\nn\\
&=&0 \quad (\mod\,2),\qquad \forall\ 1\le k\le q,\quad m\ge 1. \lb{4.7a}\eea
Now (\ref{4.6a}) and (\ref{4.7a}) complete the proof of Claim 1.

\medskip

{\bf Claim 2.} {\it Under the assumption (ECC), there holds $\hat{i}(y_k,1)=\hat{i}(y_k) > 0$, $\forall\ 1\le k\le q$.}

\medskip

At first we prove that $\hat{i}(y_k)\ge 0$, $\forall\ 1\le k\le q$. In fact, assume that there exists some $q_0\ge 1$ such that $\hat{i}(y_k) < 0$ for $1\le k\le q_0$ and $\hat{i}(y_k) \geq 0$ for $q_0+1\le k\le q$, then by the second identity in (\ref{2.11}) of Theorem 2.3, we have
\bea \sum_{1\le k\le q_0 \atop \hat{i}(y_k)<0}\frac{\hat{\chi}(y_k)}{\hat{i}(y_k)}=\sum_{1\le k\le q_0 \atop \hat{i}(y_k)<0}\frac{1}{\hat{i}(y_k)}=0,\lb{4.8a}\eea
where, in the first equality, $\hat{\chi}(y_k)=(-1)^{i(y_k)}=1$ by (\ref{2.12}) and Claim 1. This is a contradiction.

\smallskip

Now we prove $\hat{i}(y_k)\neq 0$, $\forall\ 1\le k\le q$. Otherwise, we assume that there exists some $q_0\ge 1$ such that $\hat{i}(y_k)= 0$ for $1\le k\le q_0$ and $\hat{i}(y_k) > 0$ for $q_0+1\le k\le q$. By (\ref{4.5a}) we have
\bea  \sum_{j=1}^{r_k}\frac{m\th^k_j}{2\pi}=\frac{m}{2}(r_k-i(y_k)-n-1),
    \quad\forall\  1\le k\le q_0,\ m\ge 1, \lb{4.9a}\eea
which implies
\bea  m(r_k-i(y_k)-n-1)-2r_k\le 2\sum_{j=1}^{r_k}\left[\frac{m\th^k_j}{2\pi}\right]\le m(r_k-i(y_k)-n-1), \ \forall\  1\le k\le q_0,\ m\ge 1, \lb{4.10a}\eea
Therefore, by (\ref{4.4a}), (\ref{4.10a}) and $r_k\le n-1$, it yields
\bea -2n\le -r_k-n-1 \le i(y_k^m) \le r_k-n-1\le -2,
    \quad\forall\  1\le k\le q_0,\ m\ge 1, \lb{4.11a}\eea
i.e., $\{i(y_k^m)\ |\ 1\le k\le q_0, m\in\N\}$ is a bounded set. So, for each $1\le k\le q_0$, it follows from Claim 1 that there exists an even integer $2l_k\in [-2n,-2]$ such that there exist infinitely many $m_k\in\N$ satisfying
\bea i(y_k^{m_k})=2l_k,\qquad \forall\ 1\le k\le q_0.\lb{4.12a}\eea

Now we can follow an argument in section 9 of \cite{Vit2}, for reader's conveniences we sketch the ideas below.

For large enough $a\in\R$, by Claim 1 and (\ref{4.12a}), all the closed characteristics $y_k^m$ for $1\le k\le q$ with period larger than $aT$ will have their index: either (i) equal to $2l_1$, or (ii) different from $2l_1+1$ and $2l_1-1$. Let $X^-(a,K)=\{x\in X\mid F_{a,K}(x)<0\}$ as defined in Section 7 of \cite{Vit2}.  For arbitrarily large enough $a,a^\prime$ with $a<a^\prime$, consider the exact sequence of the triple $(X,X^-(a',K),X^-(a,K))$ as follows
\bea \rightarrow H_{S^1, d(K)+2l_1}(X,X^-(a,K))&\stackrel{i_*}{\longrightarrow}& H_{S^1, d(K)+2l_1}(X,X^-(a^\prime,K))\nn\\
 &&\stackrel{\partial_*}{\longrightarrow} H_{S^1, d(K)+2l_1-1}(X^-(a^\prime,K),X^-(a,K))\rightarrow.\lb{4.13a}\eea

We claim the homomorphism $i_*$ in (\ref{4.13a}) is nonzero. In fact, on one hand, since for large enough $a<a'$, there is no any closed characteristic with index $2l_1-1$ and period locating between $aT$ and $a'T$, so by (7.25) of \cite{Vit2} and its proof, there holds $H_{S^1, d(K)+2l_1-1}(X^-(a^\prime,K),X^-(a,K))=0$, which, together with (\ref{4.13a}), implies that the homomorphism $i_*$ is surjective.

On the other hand, there holds $H_{S^1, d(K)+2l_1}(X,X^-(a',K))\neq 0$. In fact, since when $a$ increases, by (\ref{4.12a}) we always meet infinitely many closed characteristics of type (i) due to the existence of at least $y_1$ with $\hat{i}(y_1)=0$. Such closed characteristics $y_1^{m_1}$ either contribute to $H_{S^1, d(K)+2l_1}(X,X^-(a,K))$ or kill the homology class of dimension $d(K)+2l_1+1$. But the latter is zero, because there is no any closed characteristic with index $d(K)+2l_1+1$ by Claim 1, so infinitely many closed characteristics of type (i) must actually contribute to $H_{S^1, d(K)+2l_1}(X,X^-(a,K))$, which is then nonzero. By the arbitrariness of $a$ and $a'$, we obtain $H_{S^1, d(K)+2l_1}(X,X^-(a',K))\neq 0$.

In summary, it is proved that $i_*$ is surjective and $H_{S^1, d(K)+2l_1}(X,X^-(a',K))\neq 0$. So $i_*$ is nonzero. However, in Step 2 of the proof of Theorem 7.1 in \cite{Vit2}, it is shown that $i_*$ in (\ref{4.13a}) is zero. This contradiction completes the proof of Claim 2.

\medskip

{\bf Claim 3.} {\it Under the assumption (ECC), there holds $M_{2j}=b_{2j}=1$, $M_{2j+1}=b_{2j+1}=0$, $\forall\ j\in\N_0$, and $M_p=b_p=0$, $\forall\ p\le -1$.}

\medskip

In fact, by (\ref{2.14}) and $b_p=0,\ \forall\ p\le -1$, we have the following Morse inequalities
\bea M_p&\ge& b_p,\lb{4.15a}\\
M_p - M_{p-1} + \cdots +(-1)^{p}M_{0}+\cdots
&\ge& b_p - b_{p-1}+ \cdots + (-1)^{p}b_{0}, \quad\forall\ p\in\Z.\lb{4.16a}\eea
Note that the set $\{i(y_k^m)\ |\ 1\le k\le q, m\in\N\}$ is bounded below since the mean index $\hat{i}(y_k)>0$ for any $1\le k\le q$ by Claim 2. So the alternative sum $(-1)^p\sum_{j=-\infty}^{p}(-1)^j M_j$ in the left side of (\ref{4.16a}) is a finite sum by (\ref{2.13})-(\ref{2.14}) and Lemma 2.2.

Since Claim 1 implies $i(y_y^m)-i(y_k)\in 2\Z,\ \forall\, m\ge 1$ and $1\le k\le q$, by Lemma 2.2 we have
\bea M_{2j}=\sum_{k=1}^q\#\{m\ge 1\ |\ i(y_k^m)=2j\},\quad M_{2j+1}=0,\qquad \forall\ j\in\Z.\lb{4.17a}\eea

Thus it follows from $M_{2j+1}=b_{2j+1}=0$ and (\ref{4.16a}) with $p=2j+1$ that
\bea M_{2j} - M_{2j-1} + \cdots - M_1+M_{0}-\cdots
\le b_{2j} - b_{2j-1} + \cdots  - b_{1}+b_{0},\lb{4.18a}\eea
which, together with (\ref{4.16a}) for $p=2j$, implies that the equality in (\ref{4.18a}) holds. Therefore by an induction argument we obtain
\be M_{p}=b_{p}, \quad \forall\ p\in\Z,{\lb{4.19a}}\ee
which, together with $b_{2j}=1, b_{2j+1}=0$ for $j\in\N_0$ and $b_p=0$ for $p\le -1$, yields Claim 3.

\medskip

{\bf Claim 4.} {\it Under the assumptions (ECC), there holds $i(y_k)\ge 0$ and $i(y_k^{m+1})\ge i(y_k^m)+2$, $\forall\ 1\le k\le q$, $\forall\ m\ge 1$.}

\medskip

In fact, assume that there exist some $1\le k_0\le q$ such that $i(y_{k_0})\le -2$, then $M_{i(y_{k_0})}\ge 1$ by Lemma 2.2, which contradicts to $M_{i(y_{k_0})}=0$ in Claim 3 since $i(y_{k_0})\in 2\Z$. Thus $i(y_k)\ge 0$, $\forall\ 1\le k\le q$, $m\ge 1$. Now by (\ref{4.2a}) and (\ref{4.4a}) we have
\bea i(y_k^{m+1})-i(y_k^m)&=&(i(y_k)+n+1-r_k) + 2\sum_{j=1}^{r_k}\left[\frac{(m+1)\th^k_j}{2\pi}\right]-2\sum_{j=1}^{r_k}\left[\frac{m\th^k_j}{2\pi}\right] \nn\\
&\ge&(n+1)-r_k\quad \nn\\
&\ge &2,\qquad \forall\ 1\le k\le q,\quad m\ge 1, \lb{4.20a}\eea
where we use $i(y_k)\ge 0$ in the first inequality, and $r_k\le n-1$ in the second inequality.

\medskip

{\bf Proof of Theorem 1.2.}

\medskip

By the assumption (ECC) and Claim 2, we have $\hat{i}(y_k)=\hat{i}(y_k,1)>0, \forall\ 1\le k\le q$.
So for some fixed $\bar{m}\ge 1$, it follows from Theorem 3.3 that there exist infinitely many $(q+1)$-tuples $(N, m_1, \cdots, m_q)\in\N^{q+1}$ such that for any $1\le k\le q$, there holds
\bea
i(y_k,{2m_k-m}) &=& 2N-2-i(y_k,m),\quad 1\le m\le\bar{m}, \lb{4.21a}\\
i(y_k,{2m_k}) &=& 2N-1-C(M_k)+2\Delta_k,\lb{4.22a}\\
i(y_k,{2m_k+m})&=& 2N+i(y_k,m),\quad 1\le m\le\bar{m},\lb{4.23a}\eea
where, note that $S^+_{M_k}(1)=1,Q_k(m)=0$, $\forall\ m\ge 1$ by (\ref{4.1a})-(\ref{4.2a}).

Then by (\ref{4.21a})-(\ref{4.23a}) and Theorem 3.4, we obtain
\bea
 i(y_k^{2m_k-m}) &=& 2N-2n-2-i(y_k^m),\qquad \forall\ 1\le m\le \bar{m}, \lb{4.24a}\\
 i(y_k^{2m_k}) &=& 2N-C(M_k)+2\Delta_k-n-1,\lb{4.25a}\\
 i(y_k^{2m_k+m}) &=& 2N+i(y_k^m),\qquad \forall\ 1\le m\le \bar{m}.\lb{4.26a} \eea

Furthermore, by (\ref{4.24a})-(\ref{4.26a}) and Claim 4, we have
\bea
i(y_k^{m})&\le &i(y_k^{2m_k-1}) \leq 2N-2n-2,\quad \forall\ 1 \le m\le 2m_k-1,\ 1\le k\le q, \lb{4.27a}\\
i(y_k^{m})&\ge &i(y_k^{2m_k+1}) \geq 2N,\qquad \forall\ m\ge 2m_k+1,\ 1\le k\le q,\lb{4.28a}\\
i(y_k^{2m_k}) &=& 2N-C(M_k)+2\Delta_k-n-1\in [2N-2n,2N-2],\quad 1\le k\le q,\lb{4.29a}\eea
where we use $-(n-1)\le 2\Delta_k-C(M_k)\le n-1$ in the last estimate.

On one hand, by (\ref{4.27a})-(\ref{4.29a}) and Lemma 2.2, we know that the sum $\sum_{p=2N-2n}^{2N-2} M_p$ are exactly contributed by $y_k^{2m_k}$ with $1\le k\le q$, i.e., there holds
\bea \sum_{p=2N-2n}^{2N-2} M_p=q.\lb{4.30}\eea

On the other hand, it follows from Claim 3 that
\bea \sum_{p=2N-2n}^{2N-2} M_p=\sum_{p=2N-2n}^{2N-2} b_p=n.\lb{4.31}\eea

Therefore equalities (\ref{4.30}) and (\ref{4.31}) complete the proof of Theorem 1.2. \hfill\hb

\bibliographystyle{abbrv}

\begin{thebibliography}{100}
\bibitem[BLMR]{BLMR} H. Berestycki, J. M. Lasry, G. Mancini and B. Ruf, Existence of multiple periodic orbits
  on starshaped Hamiltonian systems. {\it Comm. Pure. Appl. Math.} 38 (1985), 253-289.
\bibitem[CGHHL]{CGHHL} D. Cristofaro-Gardiner, U. Hryniewicz, M. Hutchings and H. Liu,  Contact three-manifolds
with exactly two simple Reeb orbits, arXiv:2102.04970, to appear in {\it Geom. Topol.}
\bibitem[CGH]{CGH1} D. Cristofaro-Gardiner and M. Hutchings, From one Reeb orbit to two. {\it J. Diff. Geom.}  102 (2016), 25-36.
\bibitem[CGHP]{CGHP} D. Cristofaro-Gardiner, M. Hutchings, and D. Pomerleano, Torsion contact forms in three dimensions have two or infinitely many Reeb orbits, {\it Geom. Topol.} 23 (2019), 3601-3645.
\bibitem[DLW]{DLW} H. Duan, Y. Long and W. Wang, The enhanced common index jump theorem for symplectic paths and non-hyperbolic closed geodesics on Finsler manifolds. {\it Calc. Var. and PDEs.} 55 (2016), Art.145, 28 pp.
\bibitem[EkH]{EkH1} I. Ekeland and H. Hofer,  Convex Hamiltonian energy surfaces and their periodic trajectories,
  {\it Comm. Math. Phys.} 113 (1987), 419-469.
\bibitem[EkL]{EkL1} I. Ekeland and L. Lassoued,  Multiplicit\'e des trajectoires ferm\'ees d'un syst\'eme
hamiltonien sur une hypersurface d'energie convexe. {\it Ann. IHP. Anal. non lin\'eaire}. 4 (1987), 307-335.
\bibitem[GiG]{GiG1} V. Ginzburg and Y. Goren, Iterated index and the mean Euler characteristic.
{\it J. Topol. Anal.} 7 (2015), 453-481.
\bibitem[GHHM]{GHHM} V. Ginzburg, D. Hein, U. Hryniewicz and L. Macarini, Closed Reeb orbits on the
sphere and symplectically degenerate maxima. {\it Acta Math Vietnam}. 38 (2013), 55-78.
\bibitem[Gir]{Gir1} M. Girardi, Multiple orbits for Hamiltonian
systems on starshaped ernergy surfaces with symmetry. {\it Ann. IHP.
Analyse non lin\'eaire.} 1 (1984), 285-294.
\bibitem[HM]{HM} U. L. Hryniewicz and L. Macarini, Local contact homology and applications. {\it J. Topol. Anal.} 7 (2015), 167-238.
\bibitem[HuL]{HuL} X. Hu and Y. Long, Closed characteristics on non-degenerate star-shaped
hypersurfaces in $\R^{2n}$. {\it Sci. China Ser. A} 45 (2002), 1038-1052.
\bibitem[HWZ1]{HWZ1} H. Hofer, K. Wysocki and E. Zehnder, The dynamics on
three-dimensional strictly convex energy surfaces. {\it Ann. Math.} 148 (1998), 197-289.
\bibitem[HWZ2]{HWZ2} H. Hofer, K. Wysocki and E. Zehnder, Finite energy foliations of
tight three-spheres and Hamiltonian dynamics. {\it Ann. Math.} 157 (2003), 125-255.
\bibitem[LLo1]{LLo1} H. Liu and Y. Long, The existence of two closed characteristics on every compact
star-shaped hypersurface in ${\bf R}^4$. {\it Acta Math. Sinica, English Series} 32 (2016), 40-53.
\bibitem[LLo2]{LLo2} H. Liu and Y. Long, Resonance identities and stability of symmetric closed
characteristics on symmetric compact star-shaped hypersurfaces. {\it Calc. Var. and PDEs.} 54 (2015), 3753-3787.
\bibitem[LLo3]{LLo3} H. Liu and Y. Long, Irrationally elliptic closed characteristics on symmetric compact star-shaped hypersurfaces in $\R^4$. {\it  J. Fixed Point Theory Appl.}  19 (2017),  263-280.
\bibitem[LLW]{LLW14} H. Liu, Y. Long and W. Wang, Resonance indentities for closed
charactersitics on compact star-shaped hypersurfaces in $\R^{2n}$. {\it J. Funct. Anal.}
266 (2014), 5598-5638.
\bibitem[Lon1]{Lon2} Y. Long,  Bott formula of the Maslov-type index theory.
{\it Pacific J. Math.} 187 (1999), 113-149.
\bibitem[Lon2]{Lon3} Y. Long,  Precise iteration formulae of the Maslov-type index
theory and ellipticity of closed characteristics.  {\it Advances in Math.} 154 (2000), 76-131.
\bibitem[Lon3]{Lon4} Y. Long,  Index Theory for Symplectic Paths with Applications.
Progress in Math. 207, Birkh\"auser. 2002.
\bibitem[LoZ]{LoZ} Y. Long and C. Zhu, Closed characteristics on compact convex hypersurfaces
in $\R^{2n}$. {\it Ann. Math.} 155 (2002), 317-368.
\bibitem[Rab]{Rab1} P. Rabinowitz, Periodic solutions of Hamiltonian systems. {\it Comm. Pure. Appl. Math.} 31 (1978),
  157-184.
\bibitem[Szu]{Szu1} A. Szulkin, Morse theory and existence of periodic solutions of convex Hamiltonian systems. {\it
  Bull. Soc. Math. France.} 116 (1988), 171-197.
\bibitem[Vit1]{Vit1} C. Viterbo, Une th\'eorie de Morse pour les syst\`emes hamiltoniens \'etoil\'es. {\it C. R. Acad.
  Sci. Paris Ser. I Math.} 301 (1985), 487-489.
\bibitem[Vit2]{Vit2} C. Viterbo, Equivariant Morse theory for starshaped Hamiltonian systems.
{\it Trans. Amer. Math. Soc. } 311 (1989), 621-655.
\bibitem[Wan1]{Wan2} W. Wang, Existence of closed characteristics on compact convex hypersurfaces in $\R^{2n}$.
{\it Calc. Var. and PDEs.} 55 (2016), 1-25.
\bibitem[Wan2]{Wan3} W. Wang, Closed characteristics on compact convex hypersurfaces in $\R^8$. {\it Advances in Math.} 297 (2016), 93-148.
\bibitem[WHL]{WHL} W. Wang, X. Hu and Y. Long, Resonance identity, stability and multiplicity of closed
characteristics on compact convex hypersurfaces. {\it Duke Math. J.} 139 (2007), 411-462.
\bibitem[Wei]{Wei1} A. Weinstein, Periodic orbits for convex Hamiltonian systems. {\it Ann. Math.} 108 (1978), 507-518.

\end{thebibliography}

\end{document}